\documentstyle[12pt]{article}

\newcommand{\vs}{\vspace{.2in}}
\newcommand{\vsp}{\vspace{.1in}}

\newcommand{\noin}{\noindent}

\def\be{\bar {E}}

\def\p{{\bf P}}
\def\Q{{\bf Q}}
\def\C2{{\bf C}^2}
\def\c**{{{\bf C}} -\{ 0,1 \} }

\def\bk{\bar {k}}

\def\S{\Sigma}
\def\O{{\cal O}}
\def\X{{\cal X}}
\def\Z{{\cal Z}}
\def\qed{\hfill $\Box$}
\def\be{\begin{equation}}
\def\ee{\end{equation}}

\begin{document}

\title{ On a class of rational cuspidal plane curves }

\vs

\author{ H. Flenner and  M. Zaidenberg}
\date{}
\maketitle

\begin{abstract}
We obtain new examples and the complete list of the rational cuspidal plane
curves $C$ with at least three cusps, one of which has multiplicity ${\rm
deg}\,C - 2$. It occurs that these curves are projectively rigid. We also
discuss the general problem of projective rigidity of rational cuspidal plane
curves.
\end{abstract}

\vs

 A curve $C \subset \p^2$ is called {\it cuspidal} if all its singular points
are cusps.  By { \it a cusp} we mean a locally irreducible singular point. Here
we are interested in rational cuspidal plane curves. While there is a  variety
of such curves with one or two cusps [Y1-4; Sa; Ts], there are only very few
known examples with three or more cusps. The simplest one is the three cuspidal
Steiner quartic. In degree five, there are two rational cuspidal quintics with
three cusps and another one with four cusps (see [Na]). For a rational cuspidal
curve $C$ the inequality $d < 3m$ holds, where $d = {\rm deg}\,C$ and $m$ is
the maximal multiplicity of the singular points of $C$ [MaSa]. By Bezout's
theorem, $m \le d-2$ if $C$ has at least two cusps.

In this paper we will give new examples and the complete list of rational
cuspidal plane curves with at least three cusps and with $m = d-2$ (see Theorem
3.5 below). It contains all those mentioned above. Up to projective
equivalence, for any $d \ge 4$ there are exactly $[{d-1 \over 2}]$ curves of
this class. Therefore, they are all projectively rigid. We also discuss the
general problem of projective rigidity of rational cuspidal plane curves.

\vs

\section{On multiplicity sequences}

\noin {\bf 1.1. Definition.} Let $(C, \,P) \subset (\C2,\,P)$ be an irreducible
analytic plane curve germ, and let $$\C2 = V_0 \qquad
{\stackrel{\sigma_1}{\longleftarrow}} \qquad V_1 \qquad
{\stackrel{\sigma_2}{\longleftarrow}} \qquad \cdots \qquad
{\stackrel{\sigma_n}{\longleftarrow}} \qquad V_n$$ be the sequence of blow ups
over $P$ that yields the minimal embedded resolution of singularity of $C$ at
$P$. Thus, the complete preimage of $C$ in $V_n$ is a simple normal crossing
divisor $D = E + C_n$, where $E$ is the exceptional divisor of the whole
resolution and $C_n$ is the proper preimage of $C$ in $V_n$. Denote by $E_n$
the only $-1$-component of $E$, so that $E_n \cdot (D_{\rm red} - E_n) \ge 3$.

Let $E_i \subset V_i$ be the exceptional divisor of the blow up $\sigma_i$,
$C_i \subset V_i$ be the proper transform of $C$ at $V_i$, and let $P_{i-1} =
\sigma_i(E_i) \in E_{i-1} \cap C_{i-1}$ be the centrum of $\sigma_i$. Thus, $C
= C_0 \subset V_0$ and $P = P_0 \in C_0$.

Let $m_i$ denote the multiplicity of the point $P_i \in C_i$. The sequence
${\bar m}_P = (m_0,\,m_1,\dots,m_n)$, where $m_0\ge m_1 \ge \dots \ge m_n=1$,
is called {\it the multiplicity sequence of $(C,\,P)$}. We have $$\mu = 2\delta
= \sum\limits_{i=0}^n m_i(m_i - 1)\,,$$ where $\mu$ is the Milnor number of
$(C,\,P)$ and $\delta$ is the virtual number of double points of $C$ at $P$
[Mil]. \\

The following proposition gives a characterization of the multiplicity
sequences. \\

\noin {\bf 1.2. Proposition.} {\it The multiplicity sequence  ${\bar m}_P =
(m_0,\,m_1,\dots,m_n)$ has the following two properties:

\noin i) for each $i = 1,\dots,n$ there exists $k \ge 0$ such that $$m_{i-1} =
m_i + \dots +m_{i+k}\,,$$ where $$m_i = m_{i+1} =\dots =m_{i+k-1}\,,$$ and

\noin ii) if $$m_{n-r} >  m_{n-r+1}=\dots=m_n=1\,,$$ then $m_{n-r} = r-1$.

Conversely, if ${\bar m} =  (m_0,\,m_1,\dots,m_n)$ is a non--increasing
sequence of positive integers satisfying conditions i) and ii), then ${\bar m}
=  {\bar m}_P$ for some irreducible plane curve germ $(C,\,P)$. }\\

The proof is based on the following lemma.\\

\noin {\bf 1.3. Lemma.} {\it Let ${\bar m}_P = (m_0,\,m_1,\dots,m_n)$ be the
multiplicity sequence of an irreducible plane curve singularity $(C,\,P)$.
Denote by $E_i^{(k)}$ the proper transform of the exceptional divisor $E_i$ of
$\sigma_i$ at the surface $V_{i+k}$, so that, in particular, $E_i = E_i^{(0)}$.
Then the following hold.\\

\noin a) $E_iC_i = m_{i-1}$ and $$E_i^{(k)}C_{i+k} = {\rm max\,}\{0,\,m_{i-1} -
m_i -\dots - m_{i+k-1}\}\,, \,\,k> 0\,.$$ In particular, $ E_i^{(1)} C_{i+1} =
m_{i-1} - m_i$. \\

\noin b) If $$m_{i-1} > m_i + \dots + m_{i+k-1}\,,$$ then $$m_i = m_{i+1}
=\dots=m_{i+k-1}$$ and $$m_{i-1} \ge m_i + \dots + m_{i+k}\,.$$}

\noin {\bf Proof.} a) From the equalities $C_{i-1}^* := \sigma_i^*(C_{i-1}) =
C_i + m_{i-1}E_i\,,\,\,E_i^2 = -1$ and $C_{i-1}^* E_i = 0$ it follows that
$C_iE_i = m_{i-1}$. Assume by induction that a) holds for $k \le r-1$, where $r
\ge 1$. If $C_{i+r}E_i^{(r)} > 0$, then $ C_{i+r-1}E_i^{(r-1)} > 0$ and
$P_{i+r-1} \in C_{i+r-1} \cap E_i^{(r-1)}$. Therefore, by induction hypothesis
we have $$C_{i+r-1} E_i^{(r-1)} = m_{i-1} - m_i - \dots - m_{i+r-2} > 0\,,$$
$C_{i+r} = C_{i+r-1}^* - m_{i+r-1} E_{i+r}$ and $E_i^{(r)} \cdot E_{i+r} = 1$.
Hence, $E_i^{(r)} C_{i+r} = E_i^{(r)} C_{i+r-1}^* - m_{i+r-1} E_{i+r} E_i^{(r)}
= E_i^{(r-1)} C_{i+r-1} - m_{i+r-1} = m_{i-1} - m_i -\dots -m_{i+r-1}$. This
proves (a), and also proves that $$m_{i-1} \ge m_i +\dots+m_{i+r-1}$$ if
$$m_{i-1} > m_i +\dots +m_{i+r-2}\,,$$ which is the second assertion of (b).

To prove the first assertion of (b), note that $E_i^{(r-1)}$ is tangent to
$C_{i+r-1}$ at the point $P_{i+r-1}$ iff $ E_i^{(r-1)} C_{i+r-1} > m_{i+r-1}$.
As it was done in the proof of (a), one can easily show that the latter is
equivalent to the inequality $$E_i^{(r)} C_{i+r} = m_{i-1} - m_i -\dots -
m_{i+r-1} > 0\,,$$ and it implies in turn that $E_i^{(k)}$ is tangent to
$C_{i+k}$ for each $k=0,\dots,r-1$. Since by (a) $C_{i+k} E_{i+k} = m_{i+k-1}$,
the inequality $m_{i+k-1} > m_{i+k}$, where $1 \le k \le r-1$, would mean that
the curve $E_{i+k}$ is tangent to $C_{i+k}$ at $P_{i+k}$, which is impossible,
since it is transversal to $E_i^{(k)}$. Therefore, $m_{i+k-1} = m_{i+k}$ for
all $k=1, \dots, r-1$. \qed\\

\noin {\bf Proof of Proposition 1.2.} Let  ${\bar m}_P = (m_0,\,m_1,\dots,m_n)$
be the  multiplicity sequence of an irreducible plane curve singularity
$(C,\,P)$. Write $m_{i-1} = k_im_i + r_i$ with $0 \le r_i < m_i$. It follows
from Lemma 1.3(b) that $$m_i = m_{i+1} = \dots = m_{i+k_i -1}\,.$$ Thus, if
$r_i=0$, then the condition i) is fulfilled. If $r_i>0$, then $m_{i-1} > k_im_i
= m_i+\dots+ m_{i+k_i-1}$, so that by Lemma 1.3(b) we have $$m_{i-1} \ge k_im_i
+ m_{i+k_i}\,,$$ and whence $r_i \ge m_{i+k_i}$. But $r_i > m_{i+k_i}$ would
imply that $$m_{i-1} > m_i+\dots+m_{i+k_i}\,,$$ which in turn implies by Lemma
1.3(b) that $$m_i =\dots=m_{i+k_i} < r_i\,,$$ which is a contradiction.
Therefore, in this case $m_{i+k_i} = r_i$, and so $$m_{i-1} = m_i + \dots
+m_{i+k_i-1}+m_{i+k_i}\,,$$ where $$m_i=\dots=m_{i+k_i-1}\,.$$ The proof of
(ii) is easy, and so it is omited.\\

To prove the converse, we need the following lemma. For the moment we change
the convention and define the multiplicity sequences to be infinite, setting
$m_{\nu} = 1$ for all $\nu \ge n$. Thus, the sequence $(1,\,1,\,\dots)$ serves
as multiplicity sequence of a smooth germ. \\

\noin {\bf 1.4. Lemma.} {\it Let $(C,\,P)$ be an irreducible plane curve germ
with multiplicity sequence ${\bar m}_P = (m_0,\,m_1,\dots,m_n, \dots)$. Then
there exists a germ of a smooth curve $(\Gamma,\,P)$ through $P$ with $(\Gamma
C)_P = k$ iff $k$ satisfies the condition \\

\noin (*) $k = m_0 + m_1 + \dots + m_s\,\,\,$ for some $\,\,\,s > 0\,$ with
$\,\,m_0 = m_1 =
\dots = m_{s-1}\,.$}\\

\noin {\it Proof.} We proceed by induction on the number of $m_{\nu}$ which are
bigger than $1$. If it is equal to zero, i.e. if $(C,\,P)$ is a smooth germ,
then our statement is evidently true.

Let $(\Gamma,\,P) \subset (V_0,\,P)$ be a smooth curve germ through $P$, and
let $\Gamma' \subset V_1$ be the proper transform of $\Gamma$. Then $C^* = C_1
+ m_0 E_1$, and so $$k = (\Gamma C)_P = \Gamma' C_1 + m_0 \Gamma' E_1 = \Gamma'
C_1 + m_0\,.$$ If $\Gamma' C_1 = 0$, then we are done. If not, then by
induction hypothesis (applied to $C_1$) we have $$\Gamma' C_1 = m_1 + \dots +
m_s$$ for some $s > 0$ and $m_1 = \dots = m_{s-1}$. If $s = 1$ then this proves
the Lemma. If $s > 1$, i.e. $k = m_0 + m_1 + m_2 +\dots$, then we have to show
that $m_0 = m_1$. Denote by $\Gamma''$ the proper transform of $\Gamma'$ on
$V_2$. We have, as above, $$k - m_0 = \Gamma' C_1 = \Gamma'' C_2 + m_1\,,$$
which yields that  $\Gamma'' C_2 = k - m_0 - m_1 > 0$, i.e. $\Gamma''$ meets
$C_2$. Moreover, since $\Gamma' C_1 = m_1 + m_2 +\dots > m_1$, $\Gamma'$ is
tangent to $C_1$ at $P_1 \in C_1$, and hence $P_2 \in \Gamma''$. Since
$\Gamma'$ meets $E_1$ transversally, $\Gamma''$ does not meet the proper
transform $E_1^{(1)}$ of $E_1$ in $V_2$. This means that $\Gamma''$ and
$E_1^{(1)}$ meet $E_2$ in different points, and therefore $E_1^{(1)} C_2 = 0$.
By Lemma 1.3(a) we have $E_1^{(1)} C_2 = m_0 - m_1$; thus, $m_0 = m_1$. This
completes the proof in one direction.

Conversely,  assume that $k$ satisfies (*). Then $k - m_0$ satisfies (*) with
respect to $(C_1,\,P_1)$. If $k = m_0$, then any generic smooth curve $\Gamma$
through $P=P_0$ satisfies the condition $(\Gamma C)_P = k = m_0$. If $k - m_0 >
0$, then  by inductive hypothesis there is a smooth curve germ $\Gamma' \subset
V_1$ through $P_1$ with $\Gamma' C_1 = k-m_0$. Let $\Gamma$ be the image of
$\Gamma'$ in $V$. Then $\Gamma C = \Gamma' C_1 + m_0 \Gamma' E_1$. If $k - m_0
= m_1$, then $\Gamma'$ can be chosen generically, so transversally to $E_1$,
and thus we have $\Gamma C = k$. If $k - m_0 > m_1$, then as above $\Gamma''
C_2 = k - m_0 - m_1 > 0$ and so $\Gamma'' E_1^{(1)} = 0$, which implies that
$\Gamma' E_1 = 1$. Hence, $\Gamma C = k$ also in this case. The lemma is
proven. \qed \\

Returning to the proof of Proposition 1.2, fix a non-increasing sequence
${\bar m} = (m_0,\,m_1,\dots,m_n)$ that satisfies (i) and (ii). Note that the
sequence ${\bar m}' := (m_1,\dots,m_n)$ satisfies the same assumptions. Let
$\sigma_1\,:\,V_1 \to V_0 = \C2$ be the blow up at the point $P \in \C2$. Fix
a point $P_1 \in E_1 = \sigma_1^{-1} (P) \subset V_1$. Consider first the case
when $m_1 > 1$. We may assume by induction that there exists an irreducible
plane curve germ $(C_1,\,P_1)$ with multiplicity sequence ${\bar m}_{P_1} =
{\bar m}' = (m_1,\dots,m_n)$. Since $\bar m$ satisfies (i) and (ii), from
Lemma 1.4 it easily follows that there is an embedding $(C_1,\,P_1)
\hookrightarrow (V_1,\,P_1)$ such that $(E_1 C_1)_{P_1} = m_0$. Then
obviously $C := \sigma_1 (C_1) \subset \C2$ is a plane curve singularity
 with multiplicity sequence ${\bar m}_P = {\bar m} = (m_0,\,m_1,\dots,m_n)$.
Finally, assume that $m_1 =
1$. Choose $C_1 \subset V_1$ to be a smooth curve with $(C_1 E_1)_{P_1} =
m_0$. Then again $C := \sigma_1 (C_1) \subset \C2$ has multiplicity sequence
${\bar m}_P = {\bar m} = (m_0,\,m_1,\dots,m_n)$, as desired. This proves
Proposition 1.2. \qed  \\

\noin {\bf 1.5. Remark.} It is well known that the multiplicity sequence
carries the same information as the Puiseux characteristic sequence, i.e.
each of them can be computed in terms of the other [MaSa]. Moreover, the
multiplicity sequence determines the weighted dual graph of the embedded
resolution of the cusp and vice versa. This easily follows from the proofs
of (1.2) and (1.3), see also [EiNe] or [OZ1,2]. \\

\noin {\bf 1.6. } Let $f\,:\,\X \to S$ be a flat family of irreducible plane
curve singularities, i.e. there is a diagram
\begin{center}
\begin{picture}(1000,60)
\thicklines
\put(220,5){$S $}
\put(174,45){$\X$}
\put(262,45){$\C2 \times S$}
\put(220,45){$\hookrightarrow$}
\put(185,38){\vector(1,-1){20}}
\put(270,38){\vector(-1,-1){20}}
\put(175,22){$f$}
\put(270,22){pr}
\end{picture}
\end{center}
\noin and a subvariety $\S \subset \X$ such that $f\,|\,\S\,:\,\S \to S$ is
(set theoretically) bijective, $f\,|\,\X \setminus \S\,:\, \X \setminus \S \to
S$ is smooth and the fibre $X_s := f^{-1}(s)$ has a cusp at the point $\{x_s\}
= X_s \cap \S$. We say that the  family $f$ is {\it equisingular} if it
possesses a simultaneous resolution, i.e. there is a diagram
\begin{center}
\begin{picture}(1000,80)
\thicklines
\put(168,75){$\tilde \X$}
\put(268,75){$\cal Z$}
\put(220,75){$\hookrightarrow$}
\put(166,37){$\X$}
\put(250,37){$\C2 \times S$}
\put(220,37){$\hookrightarrow$}
\put(172,70){\vector(0,-1){20}}
\put(272,70){\vector(0,-1){20}}
\put(162,60){$\pi$}
\put(277,60){$\pi$}
\put(181,32){\vector(1,-1){19}}
\put(262,32){\vector(-1,-1){19}}
\put(220,3){$S$}
\put(172,17){$f$}
\put(269,17){pr}
\end{picture}
\end{center}
\noin where $\Z$ is smooth over $S$ and for each $s \in S$ the induced diagram
of the fibres
\begin{center}
\begin{picture}(1000,65)
\thicklines
\put(168,55){$\tilde X_s$}
\put(268,55){$ Z_s$}
\put(220,55){$\hookrightarrow$}
\put(167,17){$X_s$}
\put(267,17){$\C2$}
\put(220,17){$\hookrightarrow$}
\put(173,50){\vector(0,-1){20}}
\put(272,50){\vector(0,-1){20}}
\put(162,40){$\pi$}
\put(277,40){$\pi$}
\end{picture}
\end{center}
\noin yields an embedded resolution of $X_s$ in such a way that the weighted
dual graphs of $\pi^{-1} (X_s)$ are all the same.

Observe that if the family $f$ is equisingular, then all the cusps $(X_s,\,
x_s)$ have the same multiplicity sequence, see (1.5). Vice versa, we have the
following simple lemma, which will be useful in the next section.\\

\noin {\bf 1.7. Lemma. } {\it Let  $f\,:\,\X \to S$ be a flat family of
irreducible plane curve singularities. Assume that $S$ is normal and all the
cusps $(X_s,\,x_s)$, $s \in \S$, have the same multiplicity sequence. Then the
family $f$ is equisingular. } \\

\noin {\bf Proof.} Note that $\S$ is necessarily normal and $f\,|\,\S \,:\,\S
\to S$ is an isomorphism. Blowing up $\S$ gives a morphism $\pi_1\,:\,\Z_1 \to
\C2 \times S$ whose restriction to the fibre over $s$ yields the blowing up of
$\C2$ at $x_s$. Then the proper transform $\X_1$ of $\X$ in $\Z_1$ is the
blowing up $\pi\,|\,\X_1 \,:\,\X_1 \to \X$ along $\S$. The singular set of the
induced map $\X_1 \to S$ is a subvariety $\S_1$ mapped one--to--one onto $S$.
Repeating the procedure and using the fact that all multiplicity sequences of
the cusps  $(X_s,\,x_s)$ are the same, leads to a simultaneous resolution of
$f$ as above. \qed

\section{Computation of deformation invariants in terms of multiplicity
sequences}

\noin {\bf  2.1. On the Rigidity Problem.} Consider a minimal smooth completion
 $V$ of an open surface  $X = V \setminus D$ by a simple normal crossing (SNC
for short) divisor $D$. Let $\Theta_V\langle \, D \, \rangle$ be the
logarithmic tangent bundle. By [FZ] the groups $ H^i ( \Theta_V\langle \, D \,
\rangle)$ control the deformations of the pair $(V,\,D)$; more precisely, $ H^0
( \Theta_V\langle \, D \, \rangle)$ is the space of its infinitesimal
automorphisms, $ H^1 ( \Theta_V\langle \, D \, \rangle)$ is the space of
infinitesimal deformations and $ H^2 ( \Theta_V\langle \, D \, \rangle)$ gives
the obstructions for extending infinitesimal deformations. In [FZ, Lemma 1.3]
we proved that if $X$ is a $\bf Q$--acyclic  surface, i.e. $H_i(X; {\bf Q}) =
0,\,i > 0$, then the Euler characteristic of $\Theta_V\langle \, D \, \rangle$
is equal to $K_V(K_V + D)$. If, in addition, $X$ is of log--general type, i.e.
its log--Kodaira dimension $\bk (X)= 2$, then $h^0 ( \Theta_V\langle \, D \,
\rangle) = 0$ (indeed, by Iitaka's theorem [Ii,\,Theorem 6] the automorphism
group of a surface $X$ of log--general type is finite). We conjectured in [FZ]
that such surfaces are rigid and have unobstructed deformations, i.e. that for
them $$ h^1 ( \Theta_V\langle \, D \, \rangle) =  h^2 ( \Theta_V\langle \, D \,
\rangle) = 0\,,$$ and thus also $$\chi ( \Theta_V\langle \, D \, \rangle) =
0\,.$$ This, indeed, is true in all examples that we know [FZ].

Let now  $X = \p^2 \setminus C =  V \setminus D$, where $C$ is an irreducible
plane curve and $V \to \p^2$ is the minimal embedded resolution of
singularities of $C$, so that the total transform $D$ of $C$ in $V$ is an
SNC--divisor. In view of (1.6) and (1.7) the deformations of $(V,\,D)$
correspond to equisingular embedded deformations of the curve $C$ in $\p^2$. We
say shortly that $C$ is {\it projectively rigid} (resp. {\it (projectively)
unobstructed}) if the pair $(V,\,D)$ has no infinitesimal deformations, i.e. $
h^1 ( \Theta_V\langle \, D \, \rangle) =0$ (resp. $ h^2 ( \Theta_V\langle \, D
\, \rangle) = 0$)\footnote{as an abstract curve, such $C$ may have non--trivial
equisingular deformations, which might be obstructed.}.

Observe that $C \subset \p^2$ is projectively rigid iff the only equisingular
deformations of $C$ as a plane curve are those obtained via the action of the
automorphism group ${\rm PGL}\,(3,\,{\bf C})$ on $\p^2$. Indeed, suppose that
$C_t \subset \p^2,\,t \in T,$ is a family of deformations of $C_0 = C$ such
that all the members $C_t$ have at the corresponding singular points the same
multiplicity sequence. Then the singularities can be resolved simultaneously
at a family of surfaces $(V_t,\,D_t),\,t \in T$, see (1.6), (1.7). In view of
the rigidity, there is a local isomorphism with the trivial family $(V_0,\,D_0)
\times T$, and so by blowing down this leads to a family of projective
isomorphisms $C_t {\stackrel{\varphi_t}{\longrightarrow}}C_0$. The converse is
evidently true.

It is easily seen that if $C$ is a rational cuspidal curve, then the complement
$X = \p^2 \setminus C$ is $\bf Q$--acyclic. If, in addition, $C$ has at least
three cusps, then $X$ is also of log--general type [Wak]. Thus, the rigidity
conjecture of [FZ] says that such a curve $C$ should be projectively rigid and
unobstructed. Here we compute the deformation invariants of $X$ in terms of
multiplicity sequences of the cusps of $C$. In the next section we apply these
computations to check the above rigidity conjecture for the complements of
rational cuspidal curves considered there (see Lemma 3.3; cf. also section 4).
\\

\noin {\bf 2.2. Definition} (cf. [MaSa, FZ]). Let the notation be as in
Definition 1.1. The blowing up $\sigma_{i+1},\,i \ge 1$, of $V_i$ at the point
$P_i \in C_i$ is called {\it inner} (or {\it subdivisional}) if $P_i \in E_i
\cap E_{i-k}^{(k)}$ for some $k > 0$, and it is called {\it outer} (or {\it
sprouting}) in the opposite case. Note that $\sigma_1$ is neither inner nor
outer. Moreover, $\sigma_2$ is always outer, and so $\rho \ge 1$, where $\omega
= \omega_P$ resp. $\rho = \rho_P$ denotes the number of inner resp. outer
blowing ups. Denote also by $k = k_P$ the total number of blow ups, i.e. the
length of the multiplicity sequence ${\bar m}_P = (m_0,\,m_1,\dots,m_{k_P})$
minus one. Clearly, $\omega + \rho = k-1$.

By $\lceil a \rceil$ we denote the smallest integer $\ge a$. \\

\noin {\bf 2.3. Lemma.} $$\omega_P = \sum\limits_{i=1}^{k_P} (\lceil{m_{i-1}
\over m_i} \rceil - 1)$$

\noin {\bf Proof.} It is clear that the total number of exceptional curves
$E_i^{(j)} \subset V_{i+j}$, where $1 \le i+j < k$, passing through the centers
$P_{i+j}$ of the blow ups $\sigma_{i+j+1}$ is $2\omega + \rho$. If $m_{i-1} =
sm_i$, then by Lemma 1.3 $P_{i+j} \in  E_i^{(j)}$ for $j=0,\,1,\dots,s-1$, i.e.
exactly $s$ times, except in the case when $i = k_P$. If $m_{i-1} = sm_i + r$,
where $0 < r < m_i$, then this happens for $j = 0,\,1,\dots,s$, so $(s+1)$
times. In any case, this happens $\lceil{m_{i-1} \over m_i} \rceil$ times, with
the only exception when $i = k_P$. Therefore, $$2\omega + \rho =
\sum\limits_{i=1}^k \lceil{m_{i-1} \over m_i} \rceil - 1 = \sum\limits_{i=1}^k
(\lceil{m_{i-1} \over m_i} \rceil - 1) + (k - 1)\,.$$ Since $\omega + \rho =
k-1$, we have the desired result. \qed\\

\noin {\bf 2.4. Proposition.} {\it Let $V_0$ be a smooth compact complex
surface, $C \subset V_0$ be an irreducible cuspidal curve, and $V \to V_0$ be
the embedded resolution of singularities of $C$. Denote by $K_V$ resp.
$K_{V_0}$ the canonical divisor of $V$ resp. $V_0$, by $D$ the reduced total
preimage of $C$ at $V$, and by ${\bar m}_P =
(m_{P,\,0},\,m_{P,\,1},\dots,m_{P,\,k_P})$ the multiplicity sequence at $P \in
{\rm Sing}\,C$. Let, as before, $\omega_P$ be the number of inner blow ups over
$P$. Set $$\eta_P = \sum\limits_{i=0}^{k_P} (m_{P,\,i}-1)\,.$$ Then $$K_V(K_V +
D) = K_{V_0}(K_{V_0} + C) + \sum\limits_{P \in {\rm Sing}\,C} (\eta_P+\omega_P
-1)\,.$$}

\noin {\bf Proof.} Let $\sigma_{i+1}\,:\,V_{i+1} \to V_i$ be a step in the
resolution of singularities of $C$. Put $K_i = K_{V_i}$ and let $D_i$ be the
reduced total preimage of $C$ at $V_i$. We have $$K_{i+1} = K_i^* + E_{i+1}
\qquad {\rm and} \qquad D_i^* = \sigma_{i+1}^* (D_i) = D_{i+1} + (m_i -1)
E_{i+1} + \delta_i E_{i+1}\,,$$ where \[\delta_i = \left\{ \begin{array}{ll}
0 & \mbox{if  $\sigma_{i+1}$ is neither inner nor outer} \\
1 & \mbox { if $\sigma_{i+1}$ is outer} \\
2 & \mbox {if $\sigma_{i+1}$ is inner}
\end{array} \right. \] It follows that $$K_i(K_i+D_i) = K_{i+1}(K_i^* + D_i^*)
= K_{i+1} (K_{i+1} + D_{i+1} + (m_i + \delta_i -2)E_{i+1})$$ $$ =  K_{i+1}
(K_{i+1} + D_{i+1}) - (m_i + \delta_i -2)\,.$$ Thus, $$K_{i+1} (K_{i+1} +
D_{i+1}) = K_i(K_i+D_i) +(m_i -1) + (\delta_i -1)\,.$$ Now the desired equality
easily follows. \qed \\

\noin {\bf 2.5. Corollary.} {\it Let $C \subset \p^2$ be a plane cuspidal curve
of degree $d \ge 3$, and let $\pi\,:\,V \to \p^2$ be the embedded resolution of
singularities of $C$, $D$ be the reduced total preimage of $C$ in $V$ and $K =
K_V$ be the canonical divisor. Then \be \chi ( \Theta_V\langle \, D \, \rangle)
= K(K+D) = -3(d-3) + \sum\limits_{P \in {\rm Sing}\,C} (\eta_P+\omega_P
-1)\,.\ee }

\noin {\bf 2.6. Remark.}  In view of (2.5), in the case when $C \subset \p^2$
is a rational cuspidal curve with at least three cusps, the rigidity conjecture
mentioned in (2.1) in particular yields the identity $$\sum\limits_{P \in {\rm
Sing}\,C} (\eta_P+\omega_P -1) = 3(d-3)\,,$$ which, indeed, is true in all
examples that we know (see e.g. Lemma 3.3 below).

\section{Rational cuspidal plane curves of degree $d$ with a cusp of
multiplicity $d-2$}

\noin {\bf 3.1. Lemma}. {\it Let $C \subset \p^2$ be a rational cuspidal curve
of degree $d$ with a cusp $P \in C$ of multiplicity $m_P$ with multiplicity
sequence ${\bar m}_P = (m_{P,\,0},\dots, m_{P,\,k_P})$. Then the projection
$\pi_P\,:\,C \to \p^1$ from $P$ has at most $2(d-m-1)$ branching points.
Furthermore, if $Q_1,\dots,Q_s$ are the other cusps of $C$ with multiplicities
$m_1,\dots,m_s$, then $$\sum\limits_{j=1}^s (m_j - 1) + (m_{P,\,1} -1) \le
2(d-m-1)\,.$$}

\noin {\bf Proof.} By the Riemann--Hurwitz formula, applied to the composition
${\tilde \pi}_P \,:\,\p^1 = {\tilde C} \to \p^1$ of the normalization map
${\tilde C} \to C$ and the projection $\pi_P$, which has degree $d - m$, we
obtain that
$$2(d-m) = 2 + \sum\limits_{Q \in {\tilde C}} (v_Q - 1)\,,$$
where
$v_Q$ is the ramification index of ${\tilde \pi}_P$ at $Q$. The singular point
$Q_i$ of $C$ gives rise to a branching point with ramification index $\ge m_i$,
and after blowing up at $P \in C$ the first infinitesimal point to $P$ gives
rise to a branching point with ramification index $\ge m_{P,\,1}$. This proves
the lemma. \qed\\

Denote by $(m_a)$, where $m > 1$, the following multiplicity sequence: $$(m_a)
= ({\underbrace{m,\dots,m}_{a}},\,{\underbrace{1,\dots,1}_{m+1}})\,.$$  We
write simply $(m)$ instead of $(m_1)$ for $a = 1$. Notice that $(2_k)$ is the
multiplicity sequence of a simple plane curve singularity of type
$A_{2k}\,\,(x^2 + y^{2k+1} = 0)$;
thus, $(2)$ corresponds to an ordinary cusp $x^2 + y ^3 = 0$. \\

\noin {\bf 3.2. Lemma}. {\it Let $C \subset \p^2$ be a rational cuspidal curve
of degree $d$ with a cusp $P \in C$ of multiplicity $d-2$. Then $C$ has at most
three cusps. Assume further that $C$ has three cusps. Then they are not on a
line and have multiplicity sequences resp. $[(d-2),\,(2_a),\,(2_b)]$, where $a
+ b = d-2$. Each of these cusps has only one Puiseux characteristic pair; they
are, respectively, $(d-1, \,d-2),\,(2a+1, \,2),\,(2b+1,\,2)$.} \\

\noin {\bf Proof.} The projection $C \to \p^1$ from $P \in C$ being
$2$--sheeted, by the preceding Lemma it has at most two ramification points.
Thus, by Bezout's Theorem the multiplicities of other singular points are at
most two and there are at most two of them. Moreover, it follows from Lemma 3.1
that in the case when there are two more singular points, the multiplicity
sequence at $P$ should be $(d-2)$. Hence, the only multiplicity sequences in
the case of three cusps are  $[(d-2),\,(2_a),\,(2_b)]$. By the genus formula we
have $${d-2 \choose 2} + a + b = {d-1 \choose 2}\,,$$ and thus $a+b = d-2$.

That the three cusps do not lie on a line follows from Bezout's theorem. \qed\\

\noin {\bf 3.3. Lemma}. {\it Let $C \subset \p^2$ satisfies the assumptions of
Lemma 3.2. Then $C$ is projectively rigid and unobstructed \footnote{see (2.1)
for the definitions.}.}\\

\noin {\bf Proof.} Let $(V,\,D) \to (\p^2,\,C)$ be the minimal embedded
resolution of singularities of $C$. Then, first of all, the Euler
characteristic of the holomorphic tangent bundle $\chi = \chi(\Theta_V \langle
\, D \, \rangle)$ vanishes. This follows from (1). Indeed, if $P$ has
multiplicity sequence ${\bar m}_P=(m)$, then $$ \eta_P + \omega_P - 1 = 2m -
3\,,$$ whereas for the multiplicity sequence $(2_a)$ this quantity equals $a$.
Thus, under the assumptions of Lemma 3.2 we have $$\chi = 9 - 3d + (a+b) +
2(d-2) - 3 = 0\,.$$

Furthermore, the projection from the point $P \in C$ of multiplicity $d-2$
yields a
morphism $\pi_P\,:\,V \to \p^1$, which is a $\p^1$--ruling. Its restriction to
$D$ is $3$-sheeted. Moreover, $X = V \setminus D = \p^2 \setminus C$ is a
$\Q$--acyclic affine surface, i.e. $H_i (X; \,\Q) = 0,\,i>1$. By Proposition
6.2 from [FZ] it follows that $h^2 (\Theta_V \langle \, D \, \rangle) = 0$, and
so $C$ is  unobstructed. Since ${\bar k}\,(V \setminus D) = 2$ [Wak], due to
Theorem 6 from [Ii] we also have $h^0 (\Theta_V \langle \, D \, \rangle) = 0$.
Therefore, $h^1 (\Theta_V \langle \, D \, \rangle) = 0$, that means that
$(V,\,D)$ is a rigid pair, and hence $C$ is projectively rigid (see (2.1).
\qed\\

\noin {\bf 3.4. Lemma}. {\it Let $(C,\,0) \subset (\C2,\,0)$ be a plane curve
germ given parametrically by $$t \longmapsto (f(t),\,g(t)) = (t^m
,\,\sum\limits_{\nu = 1}^{\infty} c_{\nu} t^{\nu} )\,.$$ Then the multiplicity
sequence of $(C,\,0)$ has the form $$({\underbrace{m,\dots,m}_{r}},\dots)$$ iff
(**) $c_i = 0$ for all $i$ with $i < mr$ such that $m\, \not\vert \,i$.\\

\noin Furthermore, $(C,\,0)$ has multiplicity sequence $(2_r)$ iff $m=2$, the
first $r$ odd coefficients vanish: $c_1 = c_3 =\dots=c_{2r-1} =0$ and,
moreover, $c_{2r+1} \neq 0$.  }\\

\noin {\bf Proof.} After coordinate change of type $(f(t),\,g(t)) \longmapsto
(f(t),\,g(t) - p (f(t)))$, where $p \in {\bf C}  [z]$, we may assume that $c_m
= c_{2m} =\dots =c_{rm} = 0$. Then $$g(t) = c_st^s + {\rm
higher\,\,\,order\,\,\,terms}\,,$$ with $c_s \neq 0$ and either $s > rm$ or $m
\not\vert \,s$.

First of all, we show that if $(C,\,0)$ has multiplicity sequence
$({\underbrace{m,\dots,m}_{r}},\dots)$, then $s > mr$, which is equivalent to
(**). Let $s = \rho m + s_1$, where $0 \le s_1 < m$. If $\rho < r$, then after
blowing up $\rho$ times we obtain the parametrized curve germ
$$(f(t),\,g(t)/t^{\rho m})\,,$$ which still has multiplicity $m$. But since
$g(t)/t^{\rho m}$ has multiplicity $s - \rho m = s_1$, this contradicts the
assumption that $s_1 < m$. Thus, if $(C,\,0)$ has multiplicity sequence
$({\underbrace{m,\dots,m}_{r}},\dots)$, then the condition (**) is satisfied.
The converse is clear.

Finally, assume that $m = 2,\,c_1=c_3=\dots =c_{2r-1}=0$ and $c_{2r+1} \neq 0$.
Then after the above coordinate change we have $(f(t),\,g(t)) = (t^2,\,c_{2r+1}
t^{2r+1} + \dots)$, and so due to the above criterion $(C,\,0)$ has
multiplicity sequence $(2_r)$. Once again, the converse is clear. \qed\\

\noin {\bf 3.5. Theorem}. {\it For any $d \ge 4, \,a\ge b \ge 1$ with $a+b =
d-2$ there is a unique, up to projective equivalence, rational cuspidal curve
$C = C_{d,\,a} \subset \p^2$ of degree $d$ with three cusps with multiplicity
sequences $[(d-2),\,(2_a),\,(2_b)]$.

In appropriate coordinates this curve can be parametrized as $$C_{d,\,a} = (P :
Q : R) = (s^2(s-t)^{d-2}\,\,: \,\,t^2(s-t)^{d-2} \,\,:\,\,
s^2t^2q_{d,\,a}(s,\,t))\,,$$ where $q_{d,\,a}(s,\,t) = \sum\limits_{i=0}^{d-4}
c_i s^it^{d-4-i}$ and the polynomial ${\tilde q}_{d,\,a}(T) =
\sum\limits_{i=0}^{d-4} c_i T^i$ is defined as $${\tilde q}_{d,\,a}(T) =
{f_{d,\,a}(T^2) + T^{2a - 1} \over (1 + T)^{d-2}}\,.$$ Here $f_{d,\,a}(T)$ is a
polynomial of degree $d-3$ uniquely defined by the divibisility condition $(1 +
T)^{d-2} \,|\,(f_{d,\,a}(T^2) + T^{2a - 1})$.}\footnote{For the explicit
equations, see Proposition 3.9 below.} \\

\noin {\bf Proof.} Suppose that $C \subset \p^2$ is such a curve. Since by
Lemma
3.2 its three cusps are not at a line, up to projective transformation we may
assume that $C$ has cusps at the points $(0 : 0 : 1),\,(0 : 1 : 0),\,(1 : 0 :
0)$ with multiplicity sequences resp. $(d-2),\,(2_a),\,(2_b)$. Let $h = (P : Q
: R)\,:\,\p^1 \to C \hookrightarrow \p^2$ be the normalization of $C$, where
$(P : Q : R)$ is a triple of binary forms of degree $d$ without common zero
such that $$h(1 : 1) = (0 : 0 : 1)$$ $$h(0 : 1) = (0 : 1 : 0)$$ $$h(1 : 0) = (1
: 0 : 0)\,.$$ Since $C$ is required to have cusps of multiplicity $d-2$ at $h(1
: 1)$ and of multiplicity $2$ at $h(0 : 1)$ and at $h(1 : 0)$, up to
multiplication by constant factors we may write $$P(s,\,t) = (s-t)^{d-2}s^2$$
$$Q(s,\,t) = (s-t)^{d-2}t^2$$ $$R(s,\,t) = s^2t^2 q(s,t)\,,$$ where $$q(s,\,t)
= \sum\limits_{i=0}^{d-4} c_is^i t^{d-4-i}\,\,\,\,{\rm and}\,\,\,\,c_0 \neq
0,\,\,c_{d-4} \neq 0,\,\,q(1,\,1) \neq 0\,.$$ We will show that under our
assumptions $q$ is uniquely defined.

To impose the conditions that there is a cusp of type $(2_a)$ at the point $h(0
: 1) = (0 : 1 : 0)$ resp. of type $(2_b)$ at the point $h(1 : 0) = (1 : 0 :
0)$, we rewrite the above parametrization in appropriate affine coordinates at
the corresponding points. \\

\noin At $(0 : 1)$ we set $\xi = s/t$ and we have $${\tilde f}(\xi) = {P \over
Q} = {s^2 \over t^2} = \xi^2$$ $${\tilde g}(\xi) = {R \over Q} = {s^2 q(s,\,t)
\over (s-t)^{d-2}} = {\xi^2 {\tilde q}(\xi) \over (\xi - 1)^{d-2}} \,,$$ where
$${\tilde q} (\xi) =  \sum\limits_{i=0}^{d-4} c_i \xi^i \,.$$ By Lemma 3.4 $C$
has a cusp of type $(2_a)$ at $h(0 : 1) = (0 : 1 : 0)$ iff the odd coefficients
of $\xi^i$ of the function ${R \over \xi^2 Q} = {{\tilde q}(\xi) \over (\xi -
1)^{d-2}}$ vanish up to order $(2a - 3)$ (this imposes $(a-1)$ conditions) and
the coefficient of $\xi^{2a-1}$ does not vanish. \\

\noin At $(1 : 0)$ we set $\tau = t/s $ and we have $$ {\breve f}(\tau) = {Q
\over P} = {t^2 \over s^2} = \tau^2$$ $${\breve g}(\tau) = {R \over P} =
{\tau^2{\breve q} (\tau) \over (1-\tau)^{d-2}}\,,$$ where $${\breve q}(\tau) =
\sum\limits_{i=0}^{d-4} c_i\tau^{d-4-i}\,.$$ By Lemma 3.4 $C$ has a cusp of
type $(2_b)$ at $h(1 : 0) = (1 : 0 : 0)$ iff the odd coefficients of ${R \over
\tau^2 P} = {{\breve q} (\tau) \over (1-\tau)^{d-2}}$ vanish up to order
$(2b-3)$ (this imposes $(b-1)$ conditions) and the coefficient of $\tau^{2b+1}$
does not vanish.

Note that the coefficients ${\tilde c}_i$ of $\xi_i$ in ${\tilde g}(\xi)/\xi^2$
and those ${\breve c}_i$ of $\tau_i$ in ${\breve g}(\tau)/\tau^2$ are linear
functions in $c_0,\dots,c_{d-4}$. We must show that the system $${\tilde
c}_1 = {\tilde c}_3 = \dots ={\tilde c}_{2a-3} = 0,\,\,\,{\tilde c}_{2a-1} =
1$$ $${\breve c}_1 = \dots {\breve c}_{2b-3} = 0$$ has the unique solution.
Indeed, by symmetry then also  the coefficient ${\breve c}_{2b-1}$ is uniquely
defined and non--zero.
This follows from the fact that the associate homogeneous system $${\tilde c}_1
= {\tilde c}_3 = \dots ={\tilde c}_{2a-3} = {\tilde c}_{2a-1} = 0$$ $${\breve
c}_1 = \dots {\breve c}_{2b-3} = 0$$ has the unique solution, which corresponds
to $q \equiv 0$. Observe that it has $$(a-1) + (b-1) + 1 = d - 3$$ equations
and the same number of variables. To show the uniqueness we need the following
lemma. Its proof is easy and can be omited. \\

\noin {\bf 3.6. Lemma.} {\it Let $$h(T) = \sum_{\nu \ge 0} a_{\nu}T^{\nu} \in
{\bf C} [T]$$ and $${\tilde h}(T) = h(T) (1 + T^2 u(T^2))$$ for some power
series $u \in {\bf C} [[T]]$. Set ${\tilde h}(T) = \sum_{\nu \ge 0} {\tilde
a}_{\nu}T^{\nu}$. Then $${\tilde a}_1 = {\tilde a}_3 =\dots ={\tilde a}_{2k+1}
= 0$$ iff $$a_1 = a_3 =\dots = a_{2k+1} = 0\,.$$}

Returning to the proof of the theorem, put $n = d-4$ and $$F(T) = {\tilde
q}(T)(1+T)^{n+2} = {{\tilde q}(T) \over (1-T)^{n+2}}(1-T^2)^{n+2}$$ $$G(T) =
{\breve q}(T) (1+T)^{n+2} = {{\breve q}(T) \over (1-T)^{n+2}}
(1-T^2)^{n+2}\,.$$ By Lemma 3.6 the first $a$ (resp. $(b-1)$) odd coefficients
of $F(T)$ (resp. of $G(T)$) vanish iff the same is true for ${\tilde q}(T)
\over (1-T)^{n+2}$ (resp. for ${\breve q}(T) \over (1-T)^{n+2}$).
Note that by definition ${\breve q}(T) = {\tilde q}({1 \over T})T^n$. Thus, we
have that ${\rm deg}\,F = 2n+2$ is even and $$F({1 \over T})T^{2n+2} = {\tilde
q}( {1 \over T})T^n (1 + {1 \over T})^{n+2}T^{n+2} = {\breve q}(T)(1+T)^{n+2} =
G(T)\,.$$ Therefore, the conditions that the first $a$ odd coefficients of $F$
and the first $(b-1)$ odd coefficients of $G$ vanish are equivalent to $F$
being an even function: $F(T) = F(-T)$. Indeed, since $a+b-1 = d-3 = n+1$, the
above conditions mean that all odd coefficients of $F$ vanish. Now we use the
following elementary facts.\\

\noin {\bf 3.7. Lemma.} {\it Assume that $p \in {\bf C}[T]$ and $(1+T)^kp(T)$
is even. Then $(1-T)^k \,|\,p(T)$. }\\

\noin {\bf Proof.} By the condition we have $(1+T)^kp(T) = (1-T)^k p(-T)$, as
the product is even. Thus $(1-T)^k \,|\,p(T)$. \qed\\

{}From this lemma immediatly follows\\

\noin {\bf 3.8. Corollary.} {\it If ${\rm deg}\, p \le n$ and $(1+T)^{n+2}p(T)$
is even, then $p \equiv 0$.}\\

Being applied to $p = {\tilde q}$ and $F(T) = (1 + T)^{n+2}{\tilde q}(T)$,
Corollary 3.8 implies that ${\tilde q} \equiv 0$ and so $q \equiv 0$, i.e. the
above homogeneous system has a unique solution. This completes the proof of the
first part of Theorem 3.5.

As for the second one, we must prove the explicit presentation of ${\tilde q}
= {\tilde q}_{d,\,a}$. As above, it follows from the assumptions that the first
$(a-1)$ and the last $(b-1)$ odd coefficients of $F(T)$ vanish, while the
coefficient of $T^{2a-1}$ is non--zero. Therefore, $F(T) = f(T^2) + T^{2a-1}$
with $f$ being a polynomial of degree $d-3$. Hence $${\tilde q}(T) = {f(T^2) +
T^{2a-1} \over (1+T)^{d-2}}\,.$$ From the equality $F(T) = (1+T)^{d-2}{\tilde
q}(T)$ we have that $$F(-1) = F'(-1) = \dots = F^{(d-3)}(-1) = 0\,.$$ These
equations uniquely define the derivatives of the polynomial $f(\xi)$ at $\xi =
1$ up to order $(d-3)$, and therefore $f_{d,\,a}(\xi) = f(\xi) =
\sum\limits_{k=0}^{d-3} {a_k \over k!} (\xi-1)^k$ is determined in a unique
way. This completes the proof of Theorem 3.5. \qed \\

\noin {\bf 3.9. Proposition.} {\it a) The polynomial $f = f_{d,\,a}$ in Theorem
3.5 can be given as $$f(T) = \sum\limits_{k=0}^{d-3} {a_k \over k!} (T -
1)^k\,,$$ where $a_0 = 1, \,a_1 = a - {1 \over 2}$ and $$a_k = {1 \over
2^k}(2a-1)(2a-3)\dots(2a-(2k-1)) = a_1 (a_1 - 1)\dots (a_1 -
(k-1)),\,\,\,k=1,\dots,d-3\,,$$ i.e. it coincides with the corresponding
partial sum of the Taylor expansion at $T = 1$ of (the positive branch of) the
function $T^{a_1}$. \\

\noin b) In the affine chart $(X = x/z,\,Y = y/z)$ the curve $C_{d,\,a}$ as in
Theorem 3.5 can be given by the equation $p(X,\,Y) = 0$, where $p = p_{d,\,a}
\in {\bf Q}[X,\,Y]$ is defined as follows: $$p(X,\,Y) = {X^{2a+1}Y^{2b+1} - ((X
- Y)^{d-2} - XY{\hat f}(X,\,Y))^2 \over (X - Y)^{d-2}}\,,$$ and where ${\hat
f}(X,\,Y) = Y^{d-3}f({X \over Y})$ is the homogeneous polynomial which
corresponds to $f(T)$. }\\

\noin {\bf Proof.} We start with the proof of b). In the notation of Theorem
3.5 in the affine chart $\xi = s/t$ in ${\bf P}^1$ we have $${X \over Y} = {P
\over Q} = \xi^2$$ and $$X = {(\xi - 1)^{d-2} \over {\tilde q}(\xi)}\,,$$ where
$${\tilde q}(\xi) = {\tilde q}_{d,\,a}(\xi) = \sum\limits_{i=0}^{d-4} c_i
\xi_i$$ is as above. Thus, $$(\xi^2 - 1)^{d-2} = X{\tilde q}(\xi)(\xi -
1)^{d-2} = X(f_{d,\,a}(\xi^2) + \xi^{2a-1})$$ by the definition of ${\tilde
q}(\xi)$. Plugging here $\xi^2 = X/Y$ we obtain $$(X - Y)^{d-2} =
XY(Y^{d-3}f({X\over Y}) + \xi X^{a-1}Y^b) = XY{\hat f}(X,\,Y) + \xi X^a
Y^{b+1}\,.$$ Hence, $$\xi = {(X - Y)^{d-2} - XY{\hat f}(X,\,Y) \over
X^aY^{b+1}}\,$$ and so $$\xi^2 = {X \over Y} = {((X - Y)^{d-2} - XY{\hat
f}(X,\,Y))^2 \over X^{2a}Y^{2b+2}}\,.$$ Therefore, the curve $C_{d,\,a}$ in the
affine chart $(X,\,Y)$ satisfies the equation ${\tilde p} = 0$, where $${\tilde
p}(X,\,Y) = X^{2a+1}Y^{2b+1} - ((X - Y)^{d-2} - XY{\hat f}(X,\,Y))^2\,.$$ Since
$C_{d,\,a}$ is an irreducible curve of degree $d$, b) follows from the next
lemma.\\

\noin {\bf 3.10. Lemma.} $$(X-Y)^{d-2}\,|\,{\tilde p}(X,\,Y)\,.$$

\noin {\bf Proof.} We have $${\tilde p}(X,\,Y) \equiv \psi (X,\,Y)\,\,\,{\rm
mod}\, (X-Y)^{d-2}\,,$$ where $$\psi (X,\,Y) :=  X^{2a+1}Y^{2b+1} - X^2Y^2{\hat
f}^2 (X,\,Y)\,.$$ The polynomial $\psi$ is homogeneous of degree $2d-2$, and
thus it is enough to show that \be (X - 1)^{d-2} \,|\,\psi (X,\,1) \,,\ee or
equivalently, that $$ (X^2 - 1)^{d-2} \,|\,\psi (X^2,\,1) \,.$$ Since $\psi
(X^2,\,1)$ is an even polynomial and $(X^2 - 1)^{d-2} = (X - 1)^{d-2} (X +
1)^{d-2}$, by (3.7) it is sufficient to check that $$(X + 1)^{d-2}\,|\,\psi
(X^2,\,1)\,.$$ But $$\psi (X^2,\,1) = X^{4a+2} - X^4 {\hat f}^2 (X^2,\,1)
\equiv 0 \,\,\,{\rm mod} \,(X+1)^{d-2}\,,$$ because by definition, $${\hat f}
(X^2,\,1) \equiv -X^{2a-1} \,\,\,{\rm mod} \,(X+1)^{d-2}\,.$$ \qed \\

\noin {\bf Proof of Proposition 3.9, a).} From (2) it follows that
$$ f^2(T) - T^{2a-1} = (f(T) - T^{a_1})(f(T) + T^{a_1}) \equiv 0 \,\,\,{\rm
mod}\,(T-1)^{d-2} \,,$$
where by $T^{a_1}$ we mean those branch of the square root of $T^{2a-1}$ which
is positive at $T=1$. Since $(T-1)^{d-2}$ does not divide the second factor, we
have
$$ f(T) - T^{a_1} \equiv 0 \,\,\,{\rm mod}\,(T-1)^{d-2} \,.$$
Thus, indeed, $f(T)$ is the (d-3)-th partial sum of the Taylor series of the
function $T^{a_1} = T^{2a-1 \over 2}$ at the point $T=1$, and a) follows. This
proves the Proposition. \qed \\

\noin {\bf 3.11. Remark.} By the way, it follows that any rational cuspidal
plane curve $C$ with at least three cusps, one of which has multiplicity ${\rm
deg}\,C - 2$, can be defined over $\bf Q$. \\

\noin {\bf 3.12. Examples.} Here we present the affine equations $p_{d,a} = 0$
of the curves $C_{d,\,a}$ for $4 \le d \le 7$\footnote{ they were found with
"Maple".}.\\

\noin $d=4$ and $a=1$ (Steiner's quartic)
$$p_{4,3}(X,\,Y) = -{\frac {Y^{2}X^{2}}{4}}-\left (X-Y\right )^{2}+XY\left
(Y+X\right )
$$

\noin $d=5$ and $a=2$
$$p_{5,2}(X,\,Y) = {\frac {Y^{3}X^{2}}{64}}-{\frac {9\,Y^{2}X^{3}}{64}}-\left
(X-Y\right
)^{3}+XY\left ({\frac {3\,YX}{2}}-{\frac {Y^{2}}{4}}+{\frac {3\,X^{2}}
{4}}\right )
$$

\noin $d=6$ and $a=2$
$$p_{6,2}(X,\,Y) = {\frac {7\,Y^{3}X^{3}}{128}}-{\frac
{Y^{2}X^{4}}{256}}-{\frac {Y^{4}X^
{2}}{256}}-\left (X-Y\right )^{4}$$
$$+XY\left ({\frac {9\,Y^{2}X}{8}}-{
\frac {Y^{3}}{8}}+{\frac {9\,YX^{2}}{8}}-{\frac {X^{3}}{8}}\right )
$$

\noin $d=6$ and $a=3$
$$p_{6,3}(X,\,Y) = {\frac {3\,Y^{3}X^{3}}{128}}-{\frac
{25\,Y^{2}X^{4}}{256}}-{\frac {Y^{
4}X^{2}}{256}}-\left (X-Y\right )^{4}$$
$$+XY\left ({\frac {Y^{3}}{8}}-{
\frac {5\,Y^{2}X}{8}}+{\frac {15\,YX^{2}}{8}}+{\frac {5\,X^{3}}{8}}
\right )
$$

\noin $d=7$ and $a=3$
$$p_{7,3}(X,\,Y) = {\frac {475\,Y^{3}X^{4}}{16384}}-{\frac
{25\,Y^{2}X^{5}}{16384}}-{
\frac {75\,Y^{4}X^{3}}{16384}}+{\frac {9\,Y^{5}X^{2}}{16384}}-\left (X
-Y\right )^{5}$$
$$+XY\left ({\frac {3\,Y^{4}}{64}}-{\frac {5\,Y^{3}X}{16}}
+{\frac {45\,Y^{2}X^{2}}{32}}+{\frac {15\,YX^{3}}{16}}-{\frac {5\,X^{4
}}{64}}\right )
$$

\noin $d=7$ and $a=4$
$$p_{7,4}(X,\,Y) = {\frac {459\,Y^{3}X^{4}}{16384}}-{\frac
{1225\,Y^{2}X^{5}}{16384}}-{
\frac {155\,Y^{4}X^{3}}{16384}}+{\frac {25\,Y^{5}X^{2}}{16384}}-\left
(X-Y\right )^{5}$$
$$+XY\left ({\frac {7\,Y^{3}X}{16}}-{\frac {5\,Y^{4}}{64
}}-{\frac {35\,Y^{2}X^{2}}{32}}+{\frac {35\,YX^{3}}{16}}+{\frac {35\,X
^{4}}{64}}\right )\,.
$$\\

\noin {\bf 3.13. Remark.}  The weighted dual graph of the resolution of a cusp
with multiplicity sequence $(m)$ looks like

$$        
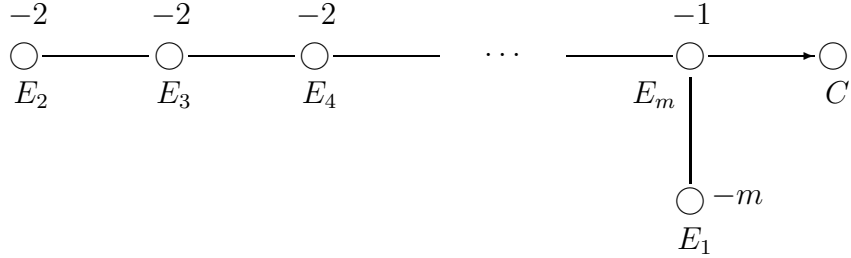
\begin{picture}(1000,90)
          \put(64,82){$-2$}
          \put(66,52){$E_2$}
          \put(70,70){\circle{10}}
          \put(77,70){\line(1,0){40}}
          \put(118,82){$-2$}
          \put(120,52){$E_3$}
          \put(125,70){\circle{10}}
          \put(132,70){\line(1,0){40}}
          \put(173,82){$-2$}
          \put(175,52){$E_4$}
          \put(180,70){\circle{10}}
          \put(187,70){\line(1,0){40}}
          \put(244,70){$\ldots$}
          \put(275,70){\line(1,0){40}}
          \put(315,82){$-1$}
          \put(300,52){$E_m$}
          \put(322,70){\circle{10}}
          \put(329,70){\vector(1,0){40}}
          \put(373,52){$C$}
          \put(376,70){\circle{10}}
          \put(322,62){\line(0,-1){40}}
          \put(322,15){\circle{10}}
          \put(330,14){$-m$}
          \put(317,-3){$E_1$}
          \end{picture}
$$
while the dual resolution graph of a cusp $(2_a) = A_{2a}$ looks like

$$      
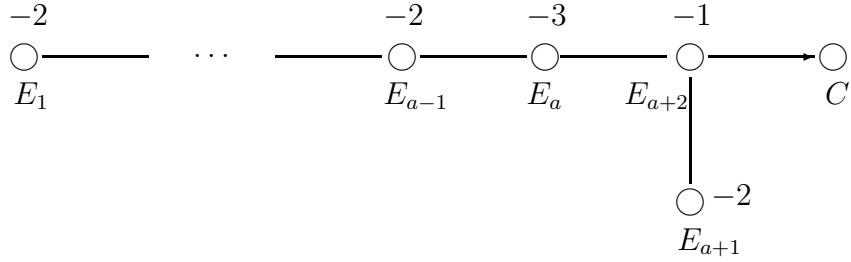
\begin{picture}(1000,90)
          \put(64,82){$-2$}
          \put(66,52){$E_1$}
          \put(70,70){\circle{10}}
          \put(77,70){\line(1,0){40}}
          \put(134,70){$\ldots$}
          \put(165,70){\line(1,0){40}}
          \put(206,82){$-2$}
          \put(206,52){$E_{a-1}$}
          \put(213,70){\circle{10}}
          \put(220,70){\line(1,0){40}}
          \put(260,82){$-3$}
          \put(260,52){$E_a$}
          \put(267,70){\circle{10}}
          \put(273,70){\line(1,0){40}}
          \put(315,82){$-1$}
          \put(297,52){$E_{a+2}$}
          \put(322,70){\circle{10}}
          \put(329,70){\vector(1,0){40}}
          \put(373,52){$C$}
          \put(376,70){\circle{10}}
          \put(322,62){\line(0,-1){40}}
          \put(322,15){\circle{10}}
          \put(330,14){$-2$}
          \put(317,-3){$E_{a+1}$}
          \end{picture}
$$\\
Therefore, the dual graph of  the total transform $D=D_{d,\,a}$ of $C_{d,\,a}$
in its minimal embedded resolution $V \to \p^2$ looks as follows:
$$  
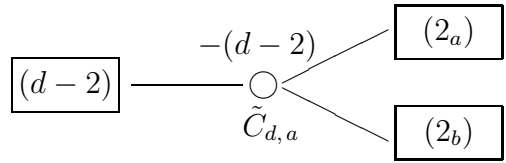
\begin{picture}(1000,60)
          \put(175,32){$-(d-2)$}
          \put(192,1){${\tilde C}_{d,\,a}$}
          \put(200,20){\circle{10}}
          \put(150,20){\line(1,0){42}}
          \put(207,21){\line(2,1){40}}
          \put(207,19){\line(2,-1){40}}
          \put(105,10){\framebox(40,20){$(d-2)$}}
          \put(250,30){\framebox(40,20){$(2_a)$}}
          \put(250,-8){\framebox(40,20){$(2_b)$}}
    \end{picture}
$$\\
where $b = d-a-2$ and boxes mean the corresponding local resolution trees, as
above. \\

\noin {\bf 3.14. Remark \footnote{This remark is due to a discussion with T.
tom Dieck, who constructed examples of cuspidal plane curves starting from
certain plane line arrangements, and with E. Artal Bartolo. We are grateful
to both of them.}.} Here we show that each curve $C_{d,a}$ can be
birationally transformed into a line. More precisely, let
$P_0, P_a, P_b$ be the cusps of
$C = C_{d,a}$ with multiplicity sequences resp. $(d-2),\, (2_a),\,(2_b)$.
Let
$l_0 = \{x = 0\}, \,l_{\infty} = \{y = 0\}$
be the lines through $P_0, P_a$, resp $P_0, P_b$, and
$l_1 = \{x - y = 0\}$ be the cuspidal tangent line to $C$ at $P_0$.
We will show that there
exist three other rational cuspidal curves $C_1, \,C_2,\,C_3$, which meet
$C$ only at the cusps of $C$, such that the curve
$T = C \cup l_0 \cup l_1 \cup l_{\infty} \cup C_1, \cup C_2 \cup C_3$
can be transformed into a configuration $T'$
of $7$ lines in $\p^2$ by means of a birational transformation $\alpha\,:\,
\p^2 \to \p^2$ which
is biregular on the complements $\p^2 \setminus T$ and $\p^2 \setminus T'$.
In fact, $\alpha$ consists of several birational
transformations composed via the following procedure.
\vsp

\noin 1) Blowing up at $P_0$, we obtain the Hirzebruch surface $\pi \,:\,
\Sigma (1) \to \p^1$
together with a two--sheeted section $C'$ (the proper preimage of $C$),
 the exceptional divisor $E$
(which is a section of $\pi$) and with
three fibres $F_0 = {l'}_0, \, F_1 = {l'}_1, \, F_{\infty} = {l'}_{\infty}$
through three points of $C'$
which we still denote resp. as $P_a,\, P_0,\, P_b$. Observe that $C'$ is smooth
at $P_0$ and by (1.3, a) $i(C',\, E;\, P_0) = d-2$.
\vsp

\noin 2) Perform $a$ resp. $b$ elementary transformations at $P_a \in C' \cap
F_0$
 resp. $P_b \in C' \cap F_{\infty}$, first blowing up at this point and then
blowing
down the proper
preimage of the fibre $F_0$ resp. $F_{\infty}$. We arrive at another
Hirzebruch surface $\Sigma (N)$ equiped with a smooth two--sheeted section
$C''$, which is
tangent to the fibres $F_0$ and $F_{\infty}$ and to the section $E'$,
where now ${E'}^2 = d-3$.
\vsp

\noin 3) Performing further $d-2$ elementary transformations at
$P_0 = E' \cap C'' \cap F_1$,
we return back at $\Sigma (1)$ with $E^2 = -1$, this time the image $C'''$
of $C''$ being a smooth two-sheeted section which does not meet $E$.
\vsp

\noin 4) Contract $E$ back to a point $P_0 \in \p^2$. Then the image
$\hat C$ of $C'''$ is a conic
in $\p^2$, and the images of the fibres $F_0,\, F_1,\, F_{\infty}$ are
resp. the lines $l_0, \,l_1,\,l_{\infty}$ through $P_0 \notin {\hat C}$,
 where $l_0,\,l_{\infty}$ are tangent to $\hat C$ resp. at the
points $P_a,\, P_b  \in {\hat C}$,
and $l_1$ is a secant line passing, say, through a point $A \in {\hat C}$.
\vsp

\noin 5) Performing the Cremona transformation with centers at the
points $A,\, P_a,\, P_b \in {\hat C}$, we obtain an arrangement $T'$ of
$7$ lines in $\p^2$ with $6$ triple points. It can be described
(in an affine chart) as a triangle together with its three medians and one
more line through the middle points of two sides. It is easily seen that
such a configuration $T'$ is projectively rigid.
\vsp

\noin The $\Q$--acyclic surface $\p^2 \setminus C$ can be
reconstructed starting from the
arrangement $T'$ by reversing the above procedure.
 In the tom Dieck-Petrie classification
[tDP, Theorem D] this line configurations is denoted as $L(4)$. \\

\noin {\bf 3.15. Remark.} E. Artal Bartolo has computed the fundamental groups
$\pi_1 (\p^2 \setminus C_{d,\,a})$. Let, as always, $a + b = d - 2$, where
$a \ge b \ge 1$. Set $2n + 1 = {\rm gcd}\,(2a+1, \,2b+1)$. Then
$\pi_1 (\p^2 \setminus C_{d,\,a}) \approx G_{d,\,n}$, where $G_{d,\,n}$ is the
group with presentation
$$G_{d,\,n} = \,<u,\,v\,|\,u(vu)^n = (vu)^n v,\,(vu)^{d-1} = v^{d-2}>\,.$$
In particular, $G_{d,\,n}$ is abelian iff $n = 0$, i.e.
${\rm gcd}\,(2a+1, \,2b+1) =1$. Furthermore, among the non--abelian groups
$G_{d,\,n}$ only $G_{4,\,1}$ and $G_{7,\,1}$ are finite. Note that, being
non--isomorphic, the curves
$C_{13,\,7}$ and $C_{13,\,10}$ have isomorphic fundamental groups of
the complements, which are both infinite non--abelian groups isomorphic to
$G_{13,\,1}$.
Evidently, there are infinitely many such pairs.

\section{Miscelleneous}

Let $C \subset {\bf P}^2$ be an irreducible plane curve, $V \to {\bf P}^2$ the
minimal embedded resolution of singularities of $C$, ${\tilde C} \subset V$ the
proper transform of $C$ and $K = K_V$ the canonical divisor of $V$. Let also $D
\subset V$ be the reduced total transform of $C$. Recall (see (2.1)) that $C$
being unobstructed simply means that  $h^2(\Theta_V\langle \, D \, \rangle) =
0$. In the next lemma we give a sufficient condition for a plane curve to be
unobstructed. \\

\noin {\bf 4.1. Lemma}. {\it Let the notation be as above.

\noin a) If $K{\tilde C} <0$, then $H^2(\Theta_V\langle \, D \, \rangle) = 0$.

\noin b) Assume that $C$ is a cuspidal curve with cusps $P_1,\dots,P_s$ having
multiplicity sequences $${\bar m}_{P_{\sigma}} = (m_{\sigma
\,1},\dots,m_{\sigma \,r_{\sigma}}, {\underbrace{1,\dots,1}_{m_{\sigma
\,r_{\sigma}} + 1}})\,,$$ where $m_{\sigma \,r_{\sigma}} \ge 2$. If $$K{\tilde
C} < \sum\limits_{\sigma = 1}^s m_{\sigma \,r_{\sigma}}\,,$$ then
$H^2(\Theta_V\langle \, D \, \rangle) = 0$.} \\

\noin {\bf Proof. } a) Fix $\omega \in H^0(\Omega^1_V \langle \, D \, \rangle
\bigotimes \omega_V)$.  Then we have ${\rm Res}_{\tilde C}(\omega) = 0 \in H^0
(\O_{\tilde C} (K{\tilde C}))$, since by assumption the degree of $\O_{\tilde
C} (K{\tilde C})$ is negative. Regarding $\omega$ as a meromorphic section in
$H^0 (\p^2,\,\Omega^1_{\p^2} \otimes \omega_{\p^2})$ it follows that $\omega$
is holomorphic outside the cusps of $C$. Therefore, $\omega$ extends to a
section in $\Omega^1_{\p^2} \otimes \omega_{\p^2}$, and hence $\omega = 0$.
Thus, $H^0(\Omega^1_V \langle \, D \, \rangle \bigotimes \omega_V) = 0$. Now
the result follows by Serre duality.

For the proof of b) consider a factorization of the embedded resolution as $$V
\to V' \to \p^2$$ such that $V' \to \p^2$ yields the minimal resolution of $C$
in the following sense:

\noin (i) The proper transform, say $C'$, of $C$ in $V'$ is smooth, and

\noin (ii) $C$ can not be resolved by fewer blowing ups.

\noin It is easily seen that $$K_{V'} C' = K_V {\tilde C} -
\sum\limits_{\sigma =
1}^s m_{\sigma \,r_{\sigma}}\,$$ (cf. the proof of (4.3, b) below). By the
above arguments, if $K_{V'} C' < 0$, then  $H^0(\Omega^1_{V'} \langle \, D'
 \, \rangle \bigotimes
\omega_{V'}) = 0$, where $D'$ is the reduced total transform of $C$ in $V'$.
Hence also $ H^0(\Omega^1_V \langle \, D \, \rangle \bigotimes \omega_V) = 0$.
\qed\\

\noin {\bf 4.2. Corollary.} {\it With the notation as in (4.1, b), assume that
$C$ is a rational cuspidal curve with ${\bar k} (\p^2 \setminus C) = 2$.
If \be \sum\limits_{\sigma = 1}^s \sum\limits_{j = 1}^{r_{\sigma}} m_{\sigma
\,j} < 3d \,,\ee then} $$\chi(\Theta_V\langle \, D \, \rangle) = K(K+D) = - h^1
(\Theta_V\langle \, D
\, \rangle) \le 0\,.$$
{\bf Proof.} From Lemma 4.3,a) below it follows that $${\tilde C}^2
+ \sum\limits_{\sigma = 1}^s m_{\sigma \,r_{\sigma}} = 3d - 2 -
\sum\limits_{\sigma = 1}^s \sum\limits_{j = 1}^{r_{\sigma}} m_{\sigma \,j}\,.
$$ Therefore, (3) is equivalent to the inequality $${\tilde C}^2 +
\sum\limits_{\sigma = 1}^s m_{\sigma \,r_{\sigma}} \ge - 1\,.$$
Thus, we have $$K{\tilde C} = -{\tilde C}^2 - 2 < \sum\limits_{\sigma =
1}^s m_{\sigma \,r_{\sigma}} \,,$$ and hence by (4.1, b)
$h^2(\Theta_V\langle \, D \, \rangle) = 0$. Since ${\bar k} (\p^2 \setminus C)
= 2$, then also $h^0(\Theta_V\langle \, D \, \rangle) = 0$  (see [Ii, Theorem
6]), and the statement follows. \qed \\

Note that in our examples, i.e. for $C = C_{d,\,a}$ being as in section 3, we
have $K_V C = d-4$ (see (4.3, b)) and $ \sum_{\sigma}
 m_{\sigma \,r_{\sigma}} = d +
2$; furthermore, $\sum\limits_{\sigma = 1}^s \sum\limits_{j = 1}^{r_{\sigma}}
m_{\sigma \,j} = 3(d-2) < 3d$. Thus, (4.1) or (4.2) gives another proof of
unobstructedness of $C_{d,\,a}$ (cf. (3.3)). \\

\noin {\bf 4.3. Lemma.} {\it Let $C \subset \p^2$ be a rational cuspidal curve,
with cusps $P_1,\dots,P_s$ having multiplicity sequences ${\bar m}_{P_{\sigma}}
= (m_{\sigma \,1},\dots,m_{\sigma \,k_{\sigma}})$. Then

\noin a) in the minimal embedded resolution $V \to \p^2$ of  singularities of
$C$ the proper transform $\tilde C$ of $C$ has selfintersection
$${\tilde C}^2
= 3d + s - 2 - \sum\limits_{i,j} m_{ij} = 3d-2 - \sum\limits_{\sigma = 1}^s
\sum\limits_{j = 1}^{r_{\sigma}} m_{\sigma \,j} - \sum\limits_{\sigma =
1}^s m_{\sigma \,r_{\sigma}}\,.$$
\noin b) Furthermore, if $K =
K_V$ is the canonical divisor, then $$K{\tilde C}= -3d - s + \sum\limits_{i,j}
m_{ij}\,.$$ }

\noin {\bf Proof.} a) Clearly, $${\tilde C}^2 = C^2 - \sum\limits_{i,j}
m_{ij}^2 + s = d^2 + s - \sum\limits_{i,j} m_{ij}^2\,.$$ The genus formula
yields $$(d-1)(d-2) =
\sum\limits_{i,j} m_{ij}(m_{ij} -1)\,.$$ Thus $$d^2 - \sum\limits_{i,j}
m_{ij}^2 = 3d - 2 - \sum\limits_{i,j} m_{ij}\,,$$ and (a) follows. \\

b) follows from (a) and the equality $K{\tilde C} + {\tilde C}^2 = -2$. An
alternative proof: we proceed by induction on the number of blow ups. First of
all, for $K = K_{{\bf P}^2}$ and $C \subset {\bf P}^2$ we have $KC = -3d$.
Furthermore, let $C \subset V$ be a curve on a surface $V$ and $K = K_V$ be the
canonical divisor, $\sigma: V' \to V$ be the blow up at a cusp of $C$ of
multiplicity $m$ and $K' = K_{V'}, \,C' \subset V'$ be the proper preimage of
$C$. We have: $$KC = K'C^* = (C' + mE)K' = C'K' + mEK' = $$ $$= C'K' + m(E(K' +
E) - E^2) = C'K' + m(-2+1) = K'C' - m\,,$$ hence $K'C' = KC + m$. This
completes the proof. \qed\\

\noin {\bf 4.4. Remark.} Let $E_P \subset V$ be the reduced exceptional divisor
of the blow ups over $P \in {\rm Sing}\,C$. Then by Lemma 2 in [MaSa] $$E_P^2 =
-\omega_P -1\,.$$ If $D = {\tilde C} + \sum\limits_{P \in {\rm Sing}\,C} E_P
\subset V$ is the reduced total transform of $C$ in $V$, then we have (cf.
[MaSa, Lemma 4]) $$D^2 = {\tilde C}^2 + 2 {\rm card}\,({\rm Sing}\,C) +
\sum\limits_{P \in {\rm Sing}\,C} E_P^2 = {\tilde C}^2 - \sum\limits_{P \in
{\rm Sing}\,C} (\omega_P -1) $$ $$ = 3d-2 - \sum\limits_{P \in {\rm Sing}\,C}
(\sum\limits_{j=0}^{k_i} m_{P,\,j} + \omega_P  - 1)\,.$$ \\

\noin {\bf 4.5. Remark.} In [OZ2, Proposition 4] the following observation is
done. \\

\noin {\it A projectively rigid rational cuspidal curve $C \subset \p^2$ cannot
have more than 9 cusps.} \\

\noin The reason is quite simple. Denote by $\kappa$  the number of cusps of
$C$. Assuming that $\kappa \ge 3$ we will have $\bk (\p^2 \setminus C) = 2$
[Wak], and therefore due to Theorem 6 in [Ii], $h^0 = 0$, where $h^i :=
h^i(\Theta_V\langle \, D \, \rangle)\,,\, i = 0,\,1,\,2$. Let $K+D = H+N$ be
the Zariski decomposition in the minimal embedded resolution $V \to \p^2$ of
singularities of $C$. It can be shown that $N^2 = \sum_{P \in {\rm Sing}\,C}
N_P^2$, where the local ingredient $N_P^2$ over a cusp $P \in  {\rm Sing}\,C$
has estimate $-N_P^2 > 1/2$. Thus,  \be {\kappa}  < 2 \sum\limits_{P \in {\rm
Sing}\,C}  (-N_P^2 ) = -2N^2 \,.\ee We also have  \be (K+D)^2 =
H^2 + N^2 \,\,\, \,\,\,{\rm and}\,\,\, \,\,\, (K+D)^2 =  K(K+D) + D(K + D) =
K(K+D) - 2 \,,\ee where [FZ, (1.3)] \be K(K + D) =
\chi(\Theta_V\langle \, D \, \rangle) =
h^2 - h^1 \,.\ee
{}From (4)--(6) and the logarithmic  Bogomolov-Miyaoka-Yau inequality
$H^2 \le
3$ [KoNaSa] we obtain $$ {\kappa}  <
-2N^2 = -2(K+D)^2 + 2H^2 \le 6  -2(K+D)^2 = 10 -  2K(K+D) = 10 - 2h^2 + 2h^1
\,.$$
Therefore, $${\kappa}  < 10 $$ as soon as $h^1 = 0$, i.e. for a projectively
rigid curve $C$. \\

Hence, once one constructs a rational cuspidal plane curve with 10 cusps or
more, we know that it is not projectively rigid. The latter means that such a
curve is a member of an equisingular \footnote{i.e. with cusps of the same
type.} family of rational cuspidal plane curves, generically pairwise
projectively non--isomorphic \footnote{i.e. non--equivalent under the action of
the automorphism group ${\rm PGL}\,(3,\,{\bf C})$ on $\p^2$.} (see (2.1)). \\

\vs
\newpage

\centerline {\bf References}

\vs

{\footnotesize

\noin [EN] D. Eisenbud, W. D. Neumann. Three-dimensional link theory and
invariants of plane curve singularities. {\it Ann. Math. Stud.}  {\bf  110},
{\it Princeton Univ. Press}, Princeton 1985

\vsp

\noin [FZ] H. Flenner, M. Zaidenberg. $\bf Q$--acyclic surfaces and their
deformations. {\it Contemporary Mathem.} {\bf 162} (1964), 143--208
\vsp

\noin [Ii] Sh. Iitaka. On logarithmic Kodaira dimension of algebraic varieties.
In: {\it Complex Analysis and Algebraic Geometry, Cambridge Univ. Press},
Cambridge e.a., 1977, 175--190
\vsp

\noin [KoNaSa] R. Kobayashi, S. Nakamura, F. Sakai. A numerical
characterization of ball quotients for normal surfaces with branch loci. {\it
Proc. Japan Acad.} {\bf  65(A)} (1989), 238--241
\vsp

\noin [MaSa] T. Matsuoka, F. Sakai. The degree of rational cuspidal curves.
{\it Math. Ann.} {\bf 285} (1989), 233--247
\vsp

\noin [Mil] J. Milnor. Singular points of complex hypersurfaces. {\it
Ann.Math.Stud.} {\bf  61}, {\it Princeton Univ. Press},  Princeton, 1968
\vsp

\noin [Na] M. Namba. Geometry of projective algebraic curves. {\it Marcel
Dekker},  N.Y. a.e., 1984
\vsp

\noin [OZ1] S.Y. Orevkov, M.G. Zaidenberg. Some estimates for plane cuspidal
curves. In: {\it Journ\'ees singuli\`eres et jacobiennes, Grenoble 26--28 mai
1993.}  Grenoble, 1994, 93--116
\vsp

\noin [OZ2] S.Y. Orevkov, M.G. Zaidenberg. On the number of singular points of
plane curves. In: {\it Algebraic Geometry. Proc. Conf., Saintama Univ., March
15--17, 1995}, 22p. (to appear)
\vsp

\noin [Sa] F. Sakai. Singularities of plane curves. {\it Preprint} (1990),
1-10
\vsp

\noin [tDP] T. tom Dieck, T. Petrie. Homology planes: An announcement
and survey. In: {\it Topological methods in algebraic transformation
groups, Progress in Mathem.} {\bf 80}, {\it Birkha\"user}, Boston e.a.,
1989, 27--48
\vsp

\noin [Ts] S. Tsunoda. The structure of open algebraic surfaces and its
application to plane curves. {\it Proc. Japan Acad.} {\bf 57(A)} (1981),
230--232
\vsp

\noin [Wak] I. Wakabayashi. On the logarithmic Kodaira dimension of the
complement of a curve in $\p^2$. {\it Proc. Japan Acad.} {\bf 54(A)} (1978),
157--162
\vsp

\noin [Wal] R. J. Walker.  Algebraic curves. {\it Princeton Univ. Press},
Princeton, 1950
\vsp

\noin [Y1] H. Yoshihara. On plane rational curves. {\it Proc. Japan Acad.} {\bf
55(A)} (1979), 152--155
\vsp

\noin [Y2] H. Yoshihara. Rational curve with one cusp. I {\it Proc. Amer. Math.
Soc.} {\bf 89} (1983), 24--26; II {\it ibid.}  {\bf 100} (1987), 405--406
\vsp

\noin [Y3] H. Yoshihara. Plane curves whose singular points are cusps. {\it
Proc. Amer. Math. Soc.} {\bf 103} (1988), 737--740

\vsp

\noin [Y4] H. Yoshihara. Plane curves whose singular points are cusps and
triple coverings of $\p^2$.  {\it Manuscr. Math.} {\bf  64} (1989), 169-187

\vs

\noin Hubert Flenner\\

\noin Fakult\"at f\"ur Mathematik\\
Ruhr Universit\"at Bochum\\
Geb.\ NA 2/72\\
Universit\"atsstr.\ 150\\
44180 BOCHUM, Germany\\

\noin e-mail: Hubert.Flenner@rz.ruhr-uni-bochum.de

\vs

\noin Mikhail Zaidenberg\\

\noin Universit\'{e} Grenoble I \\
Laboratoire de Math\'ematiques associ\'e au CNRS\\
BP 74\\
38402 St. Martin d'H\`{e}res--c\'edex, France\\

\noin e-mail: zaidenbe@fourier.grenet.fr}

\end{document}